\documentclass[10pt, a4paper]{article}

\usepackage{array}

\newcolumntype{"}{@{\hskip\tabcolsep\vrule width 1pt\hskip\tabcolsep}}

\makeatother

\usepackage{epsfig}
\usepackage{subfigure}
\usepackage{amssymb}
\usepackage{amsthm}
\usepackage{mathrsfs}
\usepackage{hhline}
\usepackage{longtable}
\usepackage{amsmath}
\usepackage{float}
\usepackage{multicol}
\usepackage{multirow}
\usepackage{cite}
\usepackage{enumerate}
\usepackage{xcolor}
\usepackage{caption}
\usepackage{tabu}

\def\ms{\medskip}
\def\nt{\noindent}

\definecolor{vividviolet}{rgb}{0.62, 0.0, 1.0}

\def\rsq{\hspace*{\fill}$\Box$\medskip}

\def\Z{\mathbb Z}
\newtheoremstyle{de}
  {10pt}          
  {10pt}  
  {\rm}  
  {}
  {\bf}  
  {. }    
  { }    
  {}     
\theoremstyle{de}

\newtheorem{de}{Definition}[section]
\newtheorem{example}{Example}[section]
\newtheorem{rem}[de]{Remark}
\newtheorem{problem}{Problem}[section]

\newtheoremstyle{theorem}
  {10pt}          
  {10pt}  
  {\it}  
  {}
  {\bf}  
  {. }    
  { }    
  {}     
\theoremstyle{theorem}
\newtheorem{theorem}{Theorem}[section]
\newtheorem{lemma}[theorem]{Lemma}

\newtheorem{corollary}[theorem]{Corollary}

\numberwithin{equation}{section}
\def\Z{\mathbb{Z}}
\def\N{\mathbb{N}}

\graphicspath{%
    {converted_graphics/}
    {/}
}

\begingroup\makeatletter\ifx\SetFigFont\undefined%
\gdef\SetFigFont#1#2#3#4#5{%
  \reset@font\fontsize{#1}{#2pt}%
  \fontfamily{#3}\fontseries{#4}\fontshape{#5}%
  \selectfont}%
\fi\endgroup%

\usepackage[mathlines]{lineno}
\modulolinenumbers[2]
\newcommand*\patchAmsMathEnvironmentForLineno[1]{%
\expandafter\let\csname old#1\expandafter\endcsname\csname #1\endcsname  \expandafter\let\csname oldend#1\expandafter\endcsname\csname end#1\endcsname  \renewenvironment{#1}%
{\linenomath\csname old#1\endcsname}%
{\csname oldend#1\endcsname\endlinenomath}}%
\newcommand*\patchBothAmsMathEnvironmentsForLineno[1]{%
\patchAmsMathEnvironmentForLineno{#1}%
\patchAmsMathEnvironmentForLineno{#1*}}%
\AtBeginDocument{%
\patchBothAmsMathEnvironmentsForLineno{equation}%
\patchBothAmsMathEnvironmentsForLineno{align}%
\patchBothAmsMathEnvironmentsForLineno{flalign}%
\patchBothAmsMathEnvironmentsForLineno{alignat}%
\patchBothAmsMathEnvironmentsForLineno{gather}%
\patchBothAmsMathEnvironmentsForLineno{multline}%
}

\begin{document}
\baselineskip18truept
\normalsize
\begin{center}
{\mathversion{bold}\Large \bf On local antimagic total chromatic number of certain\\ one point union of graphs}

\bigskip
{\large  Gee-Choon Lau}\\

\medskip

\emph{77D, Jalan Subuh, }\\
\emph{85000, Johor, Malaysia.}\\
\emph{geeclau@yahoo.com}\\

\end{center}

\begin{abstract}
Let $G = (V,E)$ be a connected simple graph of order $p$ and size $q$. A bijection $f:V(G)\cup E(G)\to \{1,2,\ldots,p+q\}$ is called a local antimagic total labeling of $G$ if for any two adjacent vertices $u$ and $v$, we have $w(u)\ne w(v)$, where  $w(u) = f(u) + \sum_{e\in E(u)} f(e)$ and $E(u)$ is the set of incident edge(s) of $u$. The local antimagic total chromatic number, denoted $\chi_{lat}(G)$, is the minimum number of distinct weights over local antimagic total labeling of $G$.  In this paper, we provide a correct proof and exact local antimagic total chromatic number of path and spider graphs given in  [Local vertex antimagic total coloring of path graph and amalgamation of path, {\it CGANT J. Maths Appln.} {\bf 5(1)} 2024, DOI:10.25037/cgantjma.v5i1.109]. Further, we determined the local antimagic total chromatic number of spider graph with each leg of length at most 2. We also showed the existence of unicyclic and bicyclic graphs with local antimagic total chromatic number 3.\\ 


\end{abstract}

\section{Introduction}
Consider graph $G=(V,E)$ of order $p$ and size $q$. In this paper, all graphs are simple and loopless. For positive integers $a < b$, let $[a,b]=\{a,a+1,\ldots,b\}$. 
%
%
Let $f: V(G)\cup E(G) \to [1,p+q]$ be a bijective total labeling that induces a vertex labeling $w_f : V(G) \to \N$, where
$$w(u)=f(u) + \sum_{uv\in E(G)} f(uv)$$
and is called the {\it weight} of $u$ for each vertex $u \in V(G)$. We say $f$ is a {\it local antimagic total labeling} of $G$ (and $G$ is {\it local antimagic total}) if $w(u)\ne w(v)$ for each $uv\in E(G)$. Clearly, $w$ corresponds to a proper vertex coloring of $G$ if each vertex $v$ is assigned the color $w(v)$. The number
\[\min\{w(f)\;|\;f \mbox{ is a local antimagic total labeling of } G\}\]
is called the {\it local antimagic total chromatic number} of $G$, denoted $\chi_{lat}(G)$. Clearly, $\chi_{lat}(G)\ge \chi(G)$.  

\ms\nt The same concepts are known as a local vertex antimagic total labeling and local vertex antimagic total chromatic number, denoted $\chi_{lvat}$, in~\cite{Putrietal}.

\ms\nt For $1\le i\le t, a_i, d_i\ge 1$ and $n = \sum^t_{i=1} d_i \ge 3$, an $n$-leg spider, denoted $Sp(a_1^{[d_1]}, a_2^{[d_2]}, \ldots, a_t^{[d_t]})$ is a tree formed by identifying an end-vertex of $n_i$ path(s) of length $d_i$. Thus, $Sp(2^{[n]})$ is the spider graph with $n$ legs of length 2. A spider graph is also know as a one point union of at least three paths at an end vertex (see~\cite{Lau+S+S}). Note that $Sp(2^{[n]})$ is known as the amalgamation of $n$ copies of $P_3$ at a pendant vertex,  denoted $Amal(P_3,x,n)$ in~\cite{Putrietal}, where $x$ is the degree $n$ vertex of $Sp(2^{[n]})$. The following two results are obtained in~\cite{Putrietal}.

\begin{theorem} For $n\ge 5$, $\chi_{lat}(P_n) = 3$. \end{theorem}

\begin{theorem} For $n\ge 3$, $\chi_{lat}(Sp(2^{[n]})) = 3$ if $n$ is odd. Otherwise, $\chi_{lat}(Sp(2^{[n]})) = 4$.   \end{theorem}

\nt Note that both $P_n$ and $Sp(2^{[n]})$ are bipartite graphs so that $\chi_{lat}(P_n), \chi_{lat}(Sp(2^{[n]}))\ge 2$. In fact, it is given in~\cite{Lau+K+S} that 

\begin{theorem} For $n\ge 2$, $\chi_{lat}(P_n) = 2$ except that $\chi_{lat}(P_4) = 3$. \end{theorem}


\begin{example} A required local antimagic total 2-labeling for $P_8$ and $P_9$ are given below.

\begin{figure}[H]
\centerline{\epsfig{file=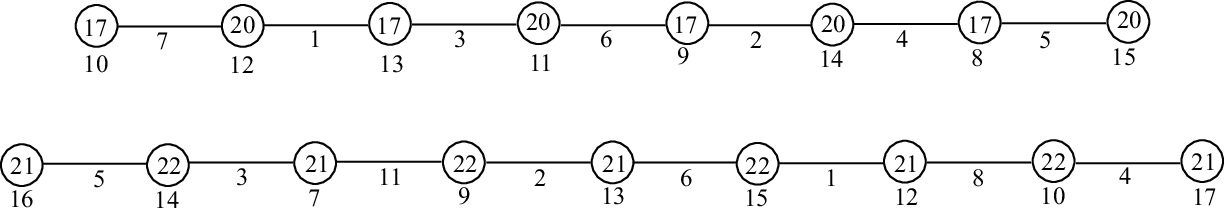, width=10cm}}
\caption{$\chi_{lat}(P_8) = \chi_{lat}(P_{9})=2$}\label{fig:1}
\end{figure}
\end{example}

\nt In this paper, we obtained a correct  proof and conclusion on $\chi_{lat}(Sp(2^{[n]}))$. Further, we determine the exact local antimagic total chromatic number of all spider graphs with each leg of length at most 2. Moreover, we show the existence of unycyclic and bicyclic graphs with local antimagic total chromatic number 3.

\section{Main Results}

We first present the correct conclusions and proof on $\chi_{lat}(Sp(2^{[n]}))$.

\begin{theorem} Suppose $n\ge 3$, $$\chi_{lat}(Sp(2^{[n]})) = \begin{cases} 2 & \mbox{ if } n\le 9\\ 3 & \mbox{ otherwise.}\end{cases} $$
\end{theorem}

\begin{proof} Let $V(Sp(2^{[n]})) = \{x, u_i, v_i\mid 1\le i\le n\}$ and $E(Sp(2^{[n]})) = \{xu_i, u_iv_i \mid 1\le i\le n\}$. Clearly, $\chi_{lat}(Sp(2^{[n]})) \ge \chi(Sp(2^{[n]})) = 2$. Suppose  $\chi_{lat}(Sp(2^{[n]})) = 2$ and $f$ is a required local antimagic total 2-coloring. We must have $w(x) = f(x) + \sum^{n}_{i=1} f(xu_i) = w(v_i) = f(u_iv_i) + f(v_i)  \ne w(u_i)$ for $1\le i\le n$. Thus, $\sum^{n+1}_{j=1} j \le w(x) = w(v_i) \le \frac{1}{n} [\sum^{4n+1}_{j=2n+2} j]$. Consequently, $(n+1)(n+2)/2\le 2n(6n+3)/2n$ so that $n^2-9n-4\le 0$. Therefore, $n\le 9$. For $3\le n\le 9$, a required local antimagic total 2-coloring of $Sp(2^{[n]})$ is shown in Figures 1 and 2 below. Thus, $\chi_{lat}(Sp(2^{[n]})) = 2$ for $3\le n\le 9$.

\begin{figure}[H]
\centerline{\epsfig{file=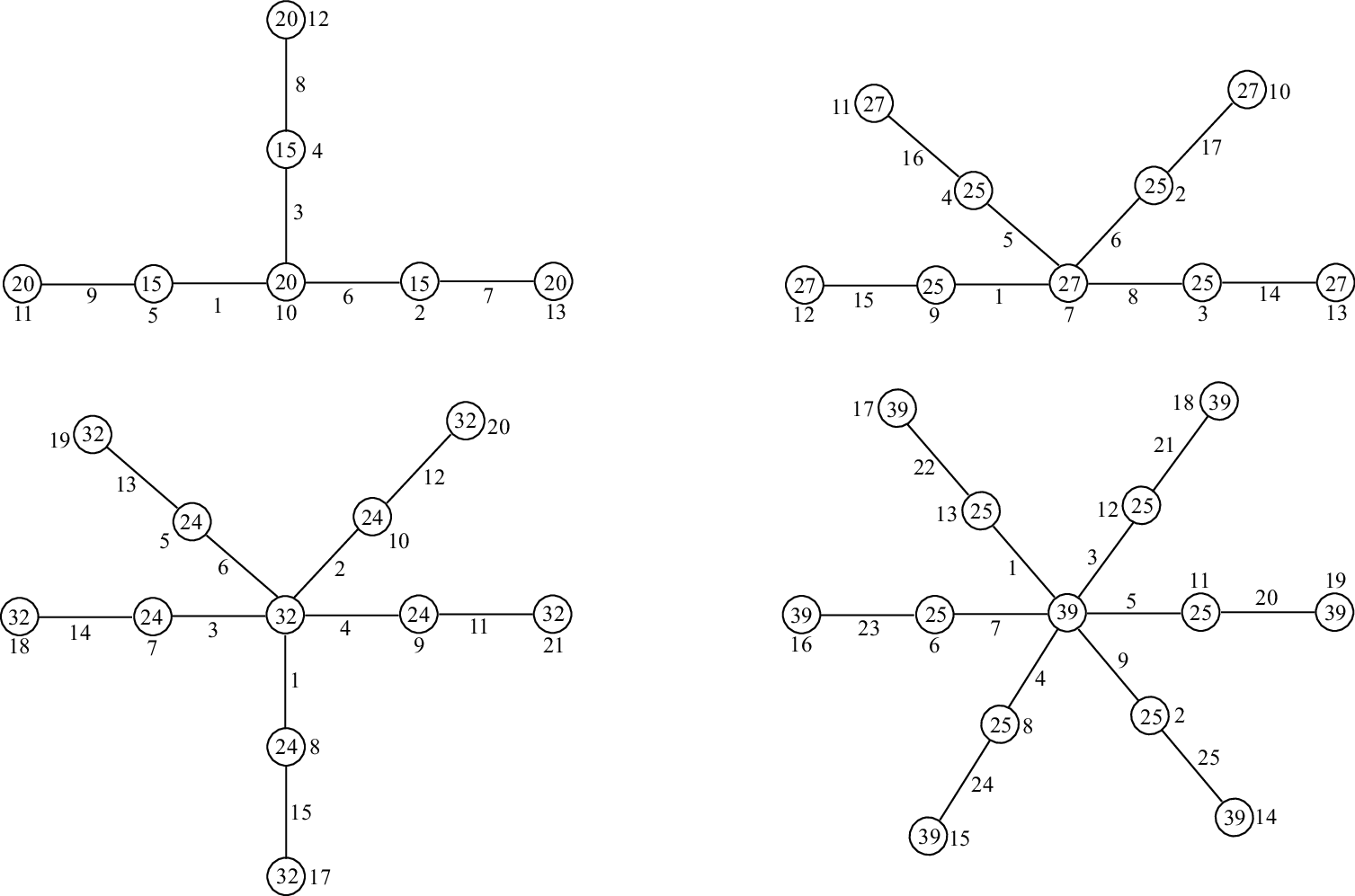, width=10cm}}
\caption{$\chi_{lat}(Sp(2^{[n]}))=2$ for $n=3,4,5,6$}\label{fig:2}
\end{figure}

\begin{figure}[H]
\centerline{\epsfig{file=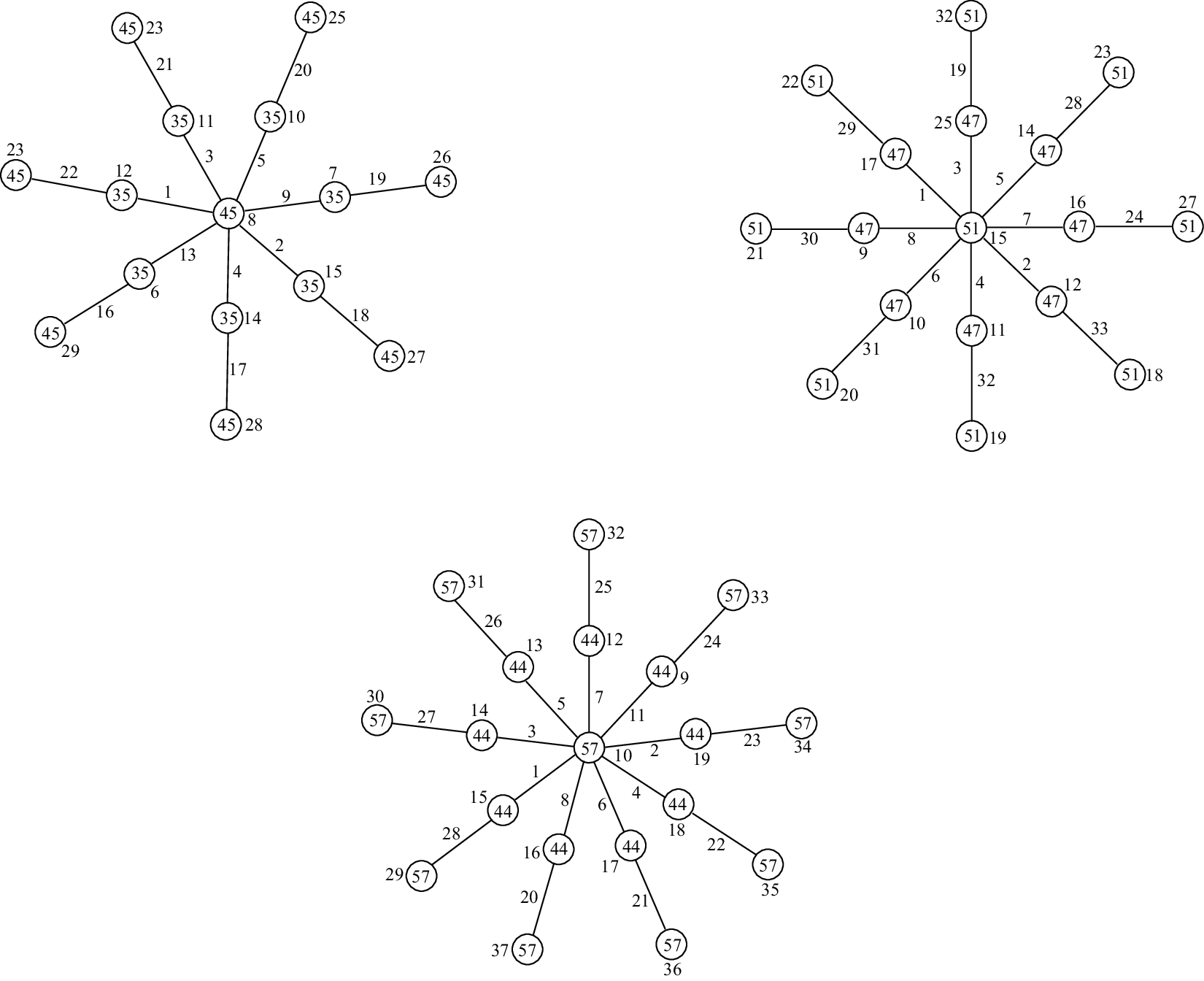, width=10cm}}
\caption{$\chi_{lat}(Sp(2^{[n]}))=2$ for $n=7,8,9$}\label{fig:3}
\end{figure}

\nt Consider $n\ge 10$. We now have $\chi_{lat}(Sp(2^{[n]}))\ge 3$. If $n$ is odd, a required local antimagic total 3-labeling is given in~\cite[Theorem 2]{Putrietal}. Suppose $n = 2k\ge 10$ is even. Define a total labeling $f : V(Sp(2^{[n]})) \cup E(Sp(2^{[n]})) \to [1, 8k+1]$ such that $f(x) = 3k+1$. For other vertices and edges, we label them as per the labeling table below.

\[\fontsize{7}{10}\selectfont
\begin{tabu}{|c|[1pt]c|c|c|c|c|c|[1pt]c|c|c|c|c|c|}\hline
i & 1 & 2 & 3 & \ldots  & k-1 & k & k+1 & k+2 & k+3 & \ldots & 2k-1 & 2k   \\\tabucline[1pt]{-}
f(xu_i) & 1 & 3 & 5 & \ldots & 2k-3 & 2k-1 & 2 & 4 & 6 & \ldots & 2k-2 & 2k  \\\hline
f(u_i) & 4k+1 & 4k & 4k-1 & \ldots & 3k+3 & 3k+2 & 3k & 3k-1 & 3K-2 & \ldots & 2k+2 & 2k+1  \\\hline
f(u_iv_i) & 5k+1 & 5k & 5k-1 & \ldots & 4k+3 & 4k+2 & 6k+1 & 6k & 6k-1 & \ldots & 5k+3 & 5k+2   \\\hline
f(v_i) & 7k+2 & 7k+3 & 7k+4 & \ldots & 8k & 8k+1 & 6k+2 & 6k+3 & 6k+4 & \ldots & 7k & 7k+1  \\\hline
\end{tabu}\]

\nt Observe that for $1\le i\le 2k$, $\{f(xu_i)\} = [1, 2k]$, $\{f(u_i)\} = [2k+1, 3k]\cup [3k+2, 4k+1]$, $\{f(u_iv_i)\} = [4k+2, 6k+1]$, and $\{f(v_i)7\} = [6k+2,8k+1]$. Thus, $f$ is a bijective total labeling of $Sp(2^{[n]})$. Observe that for each column above, the sum of the first three entries is $9k+3$ while the sum of the last two entries is $12k+3$. We now have $w(x) = 3k+1 + \sum^{2k}_{i=1} i = 2k^2 + 4k + 1$, $w(u_i) = 9k+3$, and $w(v_i) = 12k+3$. Since $w(x) \ne w(u_i) \ne w(v_i)$, $f$ is a local antimagic total 3-coloring of $Sp(2^{[n]})$ and $\chi_{lat}(Sp(2^{[n]}))\le 3$. The theorem holds. 
\end{proof}

\nt It is easy to show that $\chi_{lat}(Sp(1^{[m]})) = 2$ for $m\ge 3$. We now consider $Sp(1^{[m]}, 2^{[n]})$ for $m,n\ge 1$ and $m+n\ge 3$.

\begin{theorem}    For $(m,n) = (1,2), (2,1), (3,1), (3,2)$, $\chi_{lat}(Sp(1^{[m]}, 2^{[n]})) = 2$. For $m+n\ge 10$ or $(m,n)\in\{(6,1), (5,2), (7,1), (6,2), (5,3), (4,4), (8,1), (2,7), (3,6), (4,5), (5,4), (6,3), (7,2), (4,1), (4,3),$ $(5,1), (1,8)\}$, $\chi_{lat}(Sp(1^{[m]}, 2^{[n]})) = 3$. Otherwise, $2\le \chi_{lat}(Sp(1^{[m]}, 2^{[n]})) \le 3$ for $(m,n)\in \{(1,3), (2,2),$ $(1,4), (2,3), (4,2), (3,3), (2,4), (1,5), (3,4), (2,5), (1,6), (3,5), (2,6), (1,7)\}$.        \end{theorem}

\begin{proof} Let $G=Sp(1^{[m]}, 2^{[n]})$, $m,n\ge 1$ and $m+n\ge 3$ with $V(G) =\{x, y_i, u_j, v_j\mid 1\le i\le m, 1\le j\le n\}$ and $E(G) = \{xy_i, xu_j, u_jv_j\mid 1\le i\le m, 1\le j\le n\}$.  Consider the following two cases.

\ms\nt  {\bf Case (1).} Suppose $n = 2k +1 \ge 1$ is odd. Define a total labeling $f : V(G) \cup E(G) \to [1,2m+8k+5]$ as per the labeling table below for vertices $u_j, v_j$ and edges $xu_j, u_jv_j$, $1\le j\le 2k+1$. 

\[\fontsize{7}{10}\selectfont
\begin{tabu}{|c|[1pt]c|c|c|c|c|c|[1pt]c|c|c|c|c|c|}\hline
j & 1 & 2  & 3 & \ldots  & k & k+1 & k+2 & k+3 & k+4 & \ldots & 2k & 2k+1   \\\tabucline[1pt]{-}
f(xu_j) & 1 & 3  & 5&  \ldots & 2k-1 & 2k+1 & 2 & 4 & 6 & \ldots & 2k-2 & 2k  \\\hline
f(u_j) & 3k+2 & 3k+1 & 3k & \ldots  & 2k+3 & 2k+2 & 4k+2 & 4k+1 & 4k & \ldots & 3k+4 & 3k+3  \\\hline
f(u_jv_j) & 6k+3 & 6k+2  & 6k+1 & \ldots  & 5k+4 & 5k+3 & 5k+2 & 5k+1 & 5k &  \ldots & 4k+4 & 4k+3   \\\hline
f(v_j) & 2m+ & 2m+  & 2m+ &\ldots  & 2m+ & 2m+ & 2m+ & 2m+ & 2m+ & \ldots & 2m+ & 2m+  \\
 & 6k+4 & 6k+5 & 6k+6 &   & 7k+3 & 7k+4 & 7k+5 & 7k+6 & 7k+7 &  & 8k+3 & 8k+4 \\ \hline
\end{tabu}\]

\nt The set of entries above is $[1, 6k+3]\cup [2m+6k+4, 2m+8k+4]$. Now define $f(xy_i) = 6k+3+i$ and $f(y_i) = 2m+6k+4-i$ for $1\le i\le m$ whereas $f(x) = 2m+8k+5$ so that the corresponding set of entries is $[6k+4, 2m+6k+3]\cup\{2m+8k+5\}$. Thus, $f$ is bijective.  Observe that for each column above, the sum of the first three entries is $9k+6$ while the sum of the last two entries is $2m+12k+7$ $(= f(xy_i) + f(y_i))$. We now have $w(x) = 2m+8k+5 + (2k+1)(k+1) + \frac{1}{2}m(m+12k+7) = \frac{1}{2}m(m+12k+7) + 2m + 2k^2 + 11k + 6$, $w(y_i) = 2m + 12k+7 = w(v_j)$ and $w(u_j) = 9k+6$ for $1\le i\le m, 1\le j\le n$. Thus, $f$ is a local antimagic total 3-labeling of $G$.

\ms\nt {\bf Case (2).} Suppose $n = 2k \ge 2$ is even. Define a total labeling $f : V(G) \cup E(G) \to [1,2m+8k+1]$ as per the labeling table below for vertices $u_j, v_j$ and edges $u_jv_j$, $1\le j\le 2k$.

\[\fontsize{7}{10}\selectfont
\begin{tabu}{|c|[1pt]c|c|c|c|c|c|[1pt]c|c|c|c|c|c|}\hline
j & 1 & 2 & 3 & \ldots  & k-1 & k & k+1 & k+2 & k+3 & \ldots & 2k-1 & 2k   \\\tabucline[1pt]{-}
f(u_j) & 1 & 3 & 5 & \ldots & 2k-3 & 2k-1 & 2 & 4 & 6 & \ldots & 2k-2 & 2k  \\\hline
f(xu_j) & 2m+ & 2m+  & 2m+ &\ldots  & 2m+ & 2m+ & 2m+ & 2m+ & 2m+ & \ldots & 2m+ & 2m+  \\
  & 8k+1 & 8k & 8k-1 & \ldots & 7k+3 & 7k+2 & 7k & 7k-1 & 7k-2 & \ldots & 6k+2 & 6k+1  \\\hline
f(u_jv_j) & 3k & 3k-1 & 3k-2 & \ldots & 2k+2 & 2k+1 & 4k & 4k-1 & 4k-2 & \ldots & 3k+2 & 3k+1   \\\hline
f(v_j)   & 2m+ & 2m+  & 2m+ &\ldots  & 2m+ & 2m+ & 2m+ & 2m+ & 2m+ & \ldots & 2m+ & 2m+  \\
  & 5k+1 & 5k+2 & 5k+3 &  & 6-1k & 6k & 4k+1 & 4k+2 & 4k+3 &  & 5k-1 & 5k  \\\hline
\end{tabu}\]

\nt The set of entries above is $[1,4k]\cup [2m+4k+1, 2m+7k] \cup [2m+7k+2, 2m+8k+1]$. Now define $f(y_i) = 4k+i$ and $f(xy_i)=2m+4k+1-i$ for $1\le i\le m$ whereas $f(x) = 2m+7k+1$ so that the corresponding set of entries is $\{2m+7k+1\}\cup [4k+1,2m+4k]$. Thus, $f$ is bijective.  Observe that for each column above, the sum of the first three entries is $2m+11k+2$ while the sum of the last two entries is $2m+8k+1$ $(=f(xy_i) + f(y_i))$.  We now have $w(x) = (2m+7k+1) + \frac{1}{2}m(3m+8k+1) + k(4m+14k+2)$, $w(y_i) = 2m+8k+1 = w(v_j)$ and $w(u_j) = 2m+11k+2$ for $1\le i\le m, 1\le j\le n$. Thus, $f$ is a local antimagic total 3-labeling of $G$. 

\ms\nt Thus, $\chi_{lat}(G)\le 3$ for $m,n\ge 1$ and $m+n\ge 3$. 

\ms\nt Suppose $\chi_{lat}(G) = \chi(G) = 2$. Thus, there is a local antimagic total 2-labeling $f : V(G) \cup E(G) \to [1, 2m+4n+1]$ such that $\sum^{m+n+1}_{k=1} k \le w(x) = w(v_j) \ne w(y_i) = w(u_j) \le \frac{1}{2n}[\sum^{2m+4n+1}_{k=2m+2n+2} k]$. Therefore, $(m+n+1)(m+n+2)/2 \le 2n(4m+6n+3)/2n$ so that $(m+n+1)(m+n+2) - 2(4m+6n+3)\le 0$. Let $h = (m+n+1)(m+n+2) - 2(4m+6n+3) = (m+n)^2 - 5(m+n) - 4(n+1)$. Clearly, $h < 0$ if $m+n\le 5$. Moreover, $h > 0$ if $m+n \ge 10$. Suppose $m+n = 6$, we have $h = 4m-22 < 0$ if and only if $m\le 5$. Thus, $(m,n)\in \{(5,1), (4,2), (3,3), (2,4), (1,5)\}$. Similarly, suppose $m+n=7,8,9$, we have $h\le 0$ if and only if $(m,n) \in \{(4,3), (3,4), (2,5), (1,6), (3,5), (2,6), (1,7), (1,8)\}$.  Therefore, if $\chi_{lat}(G) = 2$, then $(m,n) \in A= \{(1,2), (2,1), (1,3), (3,1), (2,2), (1,4), (2,3), (3,2), (4,1), (5,1), (4,2), (3,3), (2,4), (1,5), (4,3), (3,4), (2,5),$ $(1,6), (3,5), (2,6), (1,7), (1,8)\}$. Consequently, $\chi_{lat}(G)  = 3$ if $m+n\ge 10$ or $(m,n) \in \{(6,1), (5,2), (7,1),$ $(6,2), (5,3), (4,4), (8,1), (2,7), (3,6), (4,5), (5,4), (6,3), (7,2)\}$.


\ms\nt We now focus on $(m,n)\in A$.  We get a local antimagic total 2-labeling for $(m,n) = (1,2), (2,1), (3,1), (3,2)$ as shown below. 

\begin{figure}[H]
\centerline{\epsfig{file=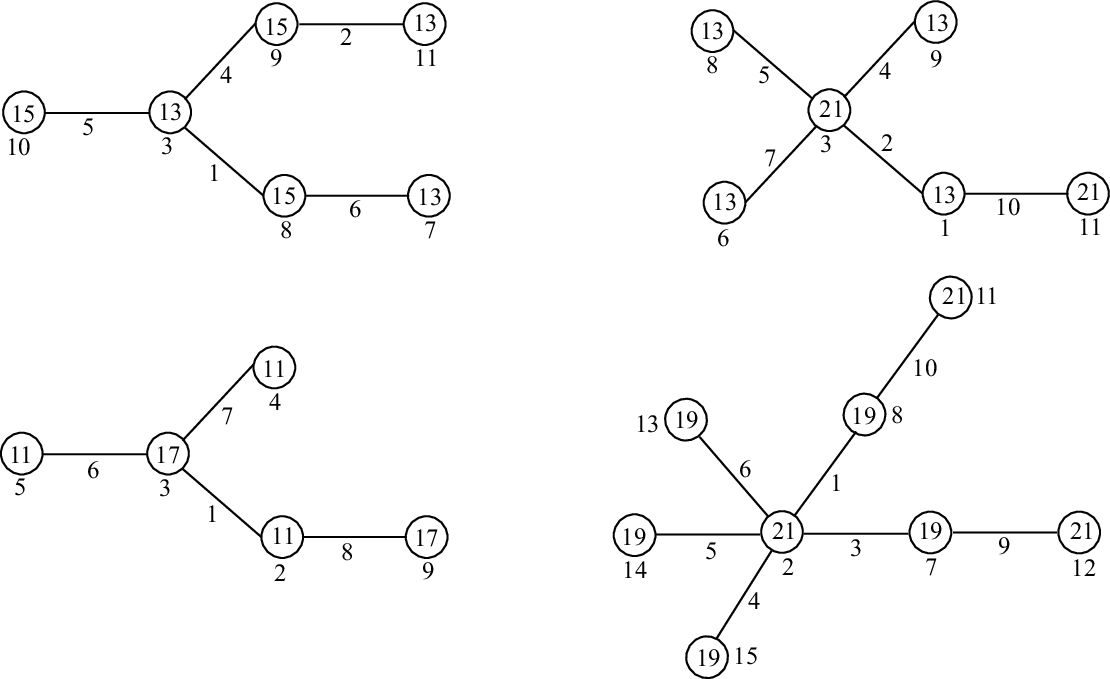, width=10cm}}
\caption{Some spider graphs with local antimagic total chromatic number 2.}\label{fig:4}
\end{figure}

\ms\nt Consider $(m,n)=(4,1)$. Suppose there is a local antimagic total 2-labeling $f$ using integers in $[1,13]$, we must have $21\le w(x) = w(v_1) \le 25$. If $w(x) = 21$, we have $\{f(u_1v_1), f(v_1)\} = \{10,11\}$ or $\{9,12\}$ or $\{8,13\}$. If $\{f(u_1v_1), f(v_1)\} = \{10,11\}$ , we must have $\{f(x), f(xy_i), f(xu_1)\mid 1\le i\le 4\} = [1,6]$ and $w(u_1)\ge 1+7+10 = 18$. However,  $\{f(xy_i), f(y_i)\mid 1\le i\le 4\} \subset [1,9]\cup\{12,13\}$ which cannot give us four pairs of integers with sum at least 18. Thus, no local antimagic 2-labeling exists. If $\{f(u_1v_1), f(v_1)\} =\{9,12\}$ or $\{8,13\}$, we also get the same conclusion. Using the same argument for $w(x) = 22,23,24,25$ also gives the same conclusion. The details are omitted. Thus, $\chi_{lat}(Sp(1^{[4]},2)) \ne 2$. By Case (1) with $n=1$ above, we can get a local antimagic total 3-labeling. Therefore, $\chi_{lat}(Sp(1^{[4]},2))=3$. Consider $(m,n) = (4,3)$. By the same argument above, we get the same conclusion. Therefore, $\chi_{lat}(Sp(1^{[4]},2^{[3]}))=3$.

\ms\nt Consider $(m,n)=(5,1)$. Suppose there is a local antimagic total 2-labeling $f$ using integers in $[1,15]$, we must have $28\le w(x) = w(v_1) \le 29$. If the former holds, we must have $\{f(x), f(xy_i), f(xu_1)\mid 1\le i\le 5\} = [1,7]$ and $\{f(u_1v_1), f(v_1)\} = \{13,15\}$. Thus, $w(y_i) = w(u_1) \ge 1 + 8 + 13 = 22$ for $1\le i\le 5$. However, $\{f(xy_i), f(y_i)\} \subset [1,12]\cup\{14\}$ which cannot give us five pairs of integers with sum at least 22. If the later holds, we must have  $\{f(x), f(xy_i), f(xu_1)\mid 1\le i\le 5\} = [1,6]\cup \{8\}$ and $\{f(u_1v_1), f(v_1)\} = \{14,15\}$. Thus, $w(y_i) = w(u_1) \ge 1 + 7+14 = 22$ for $1\le i\le 5$. However, $\{f(xy_i), f(y_i)\mid 1\le i\le 5\} \subset [1,13]$ which cannot give us five pairs of integers with sum at least 22. Thus, $\chi_{lat}(Sp(1^{[5]},2))\ne 2$. By Case (1) with $n=1$ above, we can get a local antimagic total 3-labeling. Therefore, $\chi_{lat}(Sp(1^{[5]},2)) = 3$. 

\ms\nt Consider $(m,n) = (1,8)$. Suppose there is a local antimagic total 2-labeling $f$ using integers in $[1,35]$, we must have $\{f(x), f(xy_1), f(xu_j)\mid 1\le j\le 8\} = [1,10]$ and  $\{(f(u_jv_j), f(v_j)) \mid 1\le j\le 8\} = \{(27+j, 28-j)\mid 1\le j\le 8\}$ such that $w(x) = w(v_j) = 55$ for $1\le j\le 8$. Thus, $\{f(y_1), f(u_j)\} = [11,19]$ so that $f(y_1)\le 19+10 = 29$. However, $w(y_1) = w(u_1) = \cdots = w(u_8)$ means $\sum^{8}_{j=1} w(u_j) \ge (1+\cdots + 8) + (11+\cdots + 18) + (20 +\cdots +27)$ so that $w(u_j)\ge 43$, a contradiction. Thus, $\chi_{lat}(Sp(1,2^{[8]}))\ge 3$. By Case (2) with $n=8$ above, we can get a local antimagic total 3-labeling. Therefore, $\chi_{lat}(Sp(1,2^{[8]})) = 3$.
\end{proof}

\begin{theorem} There are unicyclic and bicyclic tripartite graphs with local antimagic total chromatic number 3. \end{theorem}

\begin{proof} We shall provide a constructive approach that allows us to obtain infinitely many such graphs. Let $G$ be the one point union of an odd cycle $C_a, a\ge 3$, and a path $P_b, b\ge 2$, so that $G$ has order and size $a+b-1$. Moreover, $\chi_{lat}(G)\ge 3$. Consider $P_n, n\ge 6$ is even and the function $f : V(P_n) \cup E(P_n) \setminus\{u_1\} \to [1,2n-2]$ such that
%
%
 $$f(u_i) = \begin{cases} 2n-i & \mbox{ for odd } i\ge 3, \\ 2n-2-i & \mbox{ for even } i, i\ne n \\ 2n-2 & \mbox{ for } i=n, \end{cases}$$ and $$f(u_iu_{i+1}) = \begin{cases} \frac{i+1}{2} & \mbox{ for odd } i \ge 1, \\ \frac{n+i}{2} & \mbox{ for even } i. \end{cases}$$ Note that $f(u_1)$ is undefined. Thus, we have $w(u_1)$ is undefined, $$w(u_i) = \begin{cases} \frac{5n}{2} & \mbox{ for odd } i \ge 3, \\ \frac{5n}{2} - 2 & \mbox{ for even } i. \end{cases} $$

\nt Consider $n = a+b+1$ such that $a\ge 3$ is odd and $b\ge 2$ is even. We now begin with the path $P_{a+b+1}$ with the partial labeling $f$ defined above. Merge the vertex $u_1$ with vertex $u_{a+1}$ to get a unicyclic graph $G$ which is the one point union of a cycle $C_a$ and a path $P_b$ of order and size $a+b$ with a bijective total labeling $f : V(G) \cup E(G) \to [1, 2n-1]$ such that all the vertex weights remain unchanged except the degree 3 vertex has weight $\frac{5n}{2}-1$.   Thus, $G$ is a unicyclic tripartie graph with $\chi_{lat}(G)\le 3$. Since $\chi_{lat}(G) \ge \chi(G) = 3$, we have $\chi_{lat}(G) = 3$. 

\ms\nt Consider $P_{n}, n = a+b+1$ such that $a\ge 3$ is odd and $b\ge 4$ is even. Define a partial labeling $g : V(P_n)\cup E(P_n) \to [1,2n-2]$ such that $g(x) = f(x)$ if $x\ne u_n$ so that $g(u_1)$ and $g(u_n)$ are undefined. Now, $$w(u_i) = \begin{cases}  \frac{5n}{2} & \mbox{ for odd } i \ge 3, \\ \frac{5n}{2} - 2 & \mbox{ for even } i < n. \end{cases} $$ Merge vertices $u_1$ and $u_n$ with vertex $u_{a+1}$ to get a bicyclic graph $H$ which is the one point union of a cycle $C_a$ and another cycle $C_{b}$ with a bijective total labeling $g : V(H) \cup E(H) \to [1,2n-2]$ such that all the vertex weights remain unchanged except the degree 4 vertex has weight $3n-1$. Thus, $H$ is a bicyclic tripartite graph with $\chi_{lat}(H)\le 3$. Since $\chi_{lat}(H)\ge \chi(H) = 3$, we have $\chi_{lat}(H) = 3$. 

\ms\nt This completes the proof.     \end{proof}

\section{Conclusions and Open Problem}

In this note, we gave correct proofs and conclusions to the results given in~\cite{Putrietal}. Moreover, we determined the local antimagic total chromatic number of spider graphs with each leg of length at most 2. In~\cite[Theorem 2.4]{Lau+S+N}, the authors completely determined the local antimagic chromatic number the one poin union of cycles. We end with the following problems.

\begin{problem} Characterize spider graphs with local antimagic total chromatic number 2. \end{problem}

\begin{problem} Determine the local antimagic total chromatic number of the one point union of cycles. \end{problem}

\end{document}